\documentclass[12pt,leqno]{article}
\usepackage{amssymb,amsthm,amsmath}
\usepackage{setspace}
\usepackage{bbm}
\def\bm#1{\mathbbm{#1}}
\def\r#1{{\rm #1}}
\def\fn#1{\mathop{{\rm #1}\vphantom{\dim}}}

\makeatletter
\def\@citex[#1]#2{\if@filesw\immediate\write\@auxout{\string\citation{#2}}\fi
  \def\@citea{}\@cite{\@for\@citeb:=#2\do
    {\@citea\def\@citea{\@citesep}\@ifundefined
       {b@\@citeb}{{\bf ?}\@warning
       {Citation `\@citeb' on page \thepage \space undefined}}%
{\csname b@\@citeb\endcsname}}}{#1}}
\def\@citesep{; }
\makeatother

\newtheorem{theorem}{\indent Theorem}[section]

\newtheorem{defn}[theorem]{\indent Definition}
\newtheorem{lemma}[theorem]{\indent Lemma}
\newtheorem{prop}[theorem]{\indent Proposition}

\newtheorem{remark}{\indent Remark}[section]

\title{Groups with essential dimension one}
\author{Huah Chu$^{(1)}$, Shou-Jen Hu$^{(2)}$, Ming-chang Kang$^{(3)}$ and Jiping Zhang$^{(4)}$ \\[2mm]
\parbox{15cm}{Department of Mathematics, National Taiwan University$^{(1),(3)}$, Taipei\\
Department of Mathematics, Tamkang University$^{(2)}$, Tamshui\\
Department of Mathematics, Peking University$^{(4)}$, Beijing\\
E-mail: kang@math.ntu.edu.tw} }
\date{}

\begin{document}

\maketitle
\footnote{\hspace*{-10mm} 2000 Mathematics Subject Classification: 12E05, 12F10, 12F20, 13F30.}

\noindent Abstract. Denote by $\fn{ed}_K(G)$ the essential
dimension of $G$ over $K$. If $K$ is an algebraically closed field
with $\fn{char}K=0$, Buhler and Reichstein determine explicitly
all finite groups $G$ with $\fn{ed}_K(G)=1$ [Compositio \ Math.\
{\bf 106} (1997), Theorem 6.2]. We will prove a generalization of
this theorem when $K$ is an arbitrary field.

\newpage
\section{Introduction}

Let $K$ be an arbitrary field and $G$ be a finite groups. The
essential dimension of $G$ over $K$, denoted by $\fn{ed}_K(G)$,
was introduced by Buhler and Reichstein \cite{BR}, and was
investigated further in \cite{BF,Ka}. A related notion, the
covariant dimension, was studied in \cite{KS}. Similar notions may
be extended to the case when $G$ is an algebraic group \cite{Re}.

It is obvious that $\fn{ed}_K(G)=0$ if and only if $G=\{1\}$ the trivial group.
In \cite[Theorem 6.2]{BR} the group $G$ with $\fn{ed}_K(G)=1$ was studied.

\begin{theorem}\label{t1.1}
{\rm (Buhler and Reichstein \cite{BR})} Let $K$ be a field such
that $\fn{char}K=0$ and $K$ contains all roots of unity. If $G$ is
a nontrivial finite group, then $\fn{ed}_K(G)=1$ if and only if
$G$ is isomorphic to $\bm{Z}/n\bm{Z}$ or $D_m$ where $m$ is an odd
integer.
\end{theorem}

The purpose of this paper is to generalize the above theorem when
$K$ is an arbitrary field. The answer is the following five
theorems.

\begin{theorem}\label{t1.2}
Let $K$ be an arbitrary field. Suppose that $G$ is a nontrivial
finite group with $\fn{ed}_K(G)=1$.

\r{(1)}
If $\fn{char}K=0$,
then $G$ is isomorphic to the cyclic group $\bm{Z}/n\bm{Z}$ or
the dihedral group $D_m$ of order $2m$.

\r{(2)} If $\fn{char}K=p>0$ and $p\ne 2$, then $G$ is isomorphic
to the cyclic group $\bm{Z}/n\bm{Z}$, the dihedral group $D_m$, or
the group $G(n,p^r)$.

\r{(3)} If $\fn{char}K=2$, then $G$ is isomorphic to the cyclic
group $\bm{Z}/n\bm{Z}$, the dihedral group $D_m$, the group
$G(n,2^r)$  or the group $SL_2(\bm{F}_q)$ where $q$ is some power
of $2$.
\end{theorem}

The group $G(n,p^r)$ or $G(n,2^r)$ will be defined in Definition
\ref{d3.4} (and Formula (3.2), Lemma \ref{l3.5}). When $n=1$, the
group $G(n,p^r)$ is nothing but an elementary abelian group of
order $p^r$. We emphasize that in the definition of $G(n,p^r)$ it
is necessary that $p \nmid n$.

Because of Theorem \ref{t1.2}, it remains to find the necessary
and sufficient condition for the groups $G$ of Theorem \ref{t1.2}
to attain essential dimension one over $K$.

\begin{theorem}\label{t1.3}
Let $K$ be an arbitrary field and $G=\bm{Z}/n\bm{Z}$ be the cyclic
group of order $n$.

\r{(1)} If $\fn{char}K\nmid n$, then $\fn{ed}_K(G)=1$ if and only
if $\zeta_n+\zeta^{-1}_n\in K$ when $n$ is an odd integer, or
$\zeta_n\in K$ when $n$ is an even integer.

\r{(2)} If $\fn{char}K=p>0$ and $p\mid n$, then $\fn{ed}_K(G)=1$
if and only if $n=p$.
\end{theorem}

\begin{theorem}\label{t1.4}
Let $K$ be an arbitrary field and $G=D_n$ be the dihedral group of
order $2n$.

\r{(1)} If $\fn{char}K=0$, then $\fn{ed}_K(G)=1$ if and only if
$n$ is an odd integer and $\zeta_n+\zeta^{-1}_n\in K$.

\r{(2)} If $\fn{char}K=p>0$ and $p\ne 2$, then $\fn{ed}_K(G)=1$ if
and only if $n$ is an odd integer, $\zeta_n+\zeta^{-1}_n\in K$
when $p\nmid n$, or $n=p$ when $p\mid n$.

\r{(3)} If $\fn{char}K=2$, then $\fn{ed}_K(G)=1$ if and only if
$\zeta_n+\zeta^{-1}_n\in K$ when $n$ is an odd integer, or $|K|\ge
4$ with $n=2$ when $n$ is an even integer.
\end{theorem}

\begin{theorem}\label{t1.5}
Let $K$ be an arbitrary field with $\fn{char}K=p>0$. If $G$ is the
group $G(n,p^r)$, then $\fn{ed}_K(G)=1$ if and only if $n$ is an
odd integer, $\zeta_n\in K$ and $[K: \bm{F}_p]\ge r$.
\end{theorem}

\begin{theorem}\label{t1.6}
Let $K$ be an arbitrary field with $\fn{char}K=2$.
If $G$ is the group $SL_2(\bm{F}_q)$ where $q$ is some power of $2$,
then $\fn{ed}_K(G)=1$ if and only if $K\supset \bm{F}_q$.
\end{theorem}

As an application of the above theorems, we will prove that, when
$K$ is a field with $\fn{char}K=2$, if $K$ doesn't contain
$\bm{F}_4$, then $\fn{ed}_K(A_4)=\fn{ed}_K(A_5)=2$, while
$\fn{ed}_K(A_4)=\fn{ed}_K(A_5)=1$ if $K \supset \bm{F}_4$ (see
Proposition \ref{p7.4}). Similarly, since $\bm{Z}/4\bm{Z}$ is
contained in the symmetric group $S_4$ and $\fn{ed}_K(S_4)= 2$, we
find that $\fn{ed}_K(\bm{Z}/4\bm{Z})=2$ if $\fn{char}K \neq 2$ and
$\sqrt{-1} \notin K$; $\fn{ed}_K(\bm{Z}/4\bm{Z})=1$ if $\fn{char}K
\neq 2$ and $\sqrt{-1} \in K$; $\fn{ed}_K(\bm{Z}/4\bm{Z})=2$ if
$\fn{char}K = 2$. (This result was proved in \cite[Theorem
7.6]{BF} in the case $\fn{char}K \neq 2$ by a different method.)
It is not difficult to verify that
$\fn{ed}_{\bm{Q}}(\bm{Z}/5\bm{Z})=\fn{ed}_{\bm{Q}}(\bm{Z}/6\bm{Z})=2$
by the same way; we leave the details to the reader.

We will organize the article as follows. We recall the definition
of essential dimensions and prove some basic facts in Section 2.
In Section 3 we will define several groups which arise as
subgroups of $SL_2(K)$ where $K$ is any algebraically closed field
($\fn{char}K$ may be zero or positive). In particular, the
definition of $G(n, p^r)$ in Theorem \ref{t1.2} and Theorem
\ref{t1.5} will be given. The proof of Theorem \ref{t1.2}, Theorem
\ref{t1.3}, Theorem \ref{t1.4}, Theorem \ref{t1.5} and Theorem
\ref{t1.6} will be given in Section 4, Section 5, Section 6,
Section 7 and Section 8 respectively.

Standing notation and terminology. For emphasis, $K$ is an
arbitrary field. All the fields in this article are assumed to
contain the ground field $K$. If $E$ is a field extension of $K$,
$\fn{trdeg}_K E$ denotes the transcendence degree of $E$ over $K$.
If $F\subset K$ are fields, $[K:F]$ denotes the vector space
dimension of $K$ over $F$. The order of an element $\sigma$ in a
group $G$ is denoted by $\fn{ord}(\sigma)$. All the groups $G$ in
the sequel are nontrivial finite groups. When we talk about $G\to
GL(V)$ is a representation of a finite group $G$, it is understood
that $V$ is a finite-dimensional vector space over $K$. We will
adopt the following notations,

\begin{quotation}
$S_n$, the symmetric group,

$A_n$, the alternating group,

$PGL_2(K)$, the group isomorphic to $GL_2(K)/K^{\times}$,

$\bm{Z}/n\bm{Z}$, the cyclic group of order $n$,

$D_n$, the dihedral group of order $2n$.
\end{quotation}

We will take the convention that $\fn{char}K\nmid  n$ means either
$\fn{char}K=0$ or $\fn{char}K=p>0$ with $p\nmid n$. $\zeta_n$
denotes a primitive $n$-th root of unity; whenever we write
$\zeta_n \in K$, it is assumed tacitly that $\fn{char}K\nmid n$.
Finally $\bm{F}_q$ is the finite field consisting of $q$ elements.

\section{Preliminaries}

Throughout this paper,
$K$ is an arbitrary field.
All the fields in this paper are extension fields of $K$.

\begin{defn}\label{d2.1}\rm
Let $G$ be a finite group and $L$ be a field containing $K$.
We will call $L$ a $G$-field (over $K$) if $G$ acts on $L$ by $K$-automorphism;
$L$ is a faithful $G$-field if the group homomorphism $G\to \fn{Aut}_K(L)$ is injective.
\end{defn}

\begin{defn}\label{d2.2}\rm
Let $G$ be a finite group and $K$ be an arbitrary field. Let
$\rho:G\to GL(V)$ be a faithful finite-dimensional representation
of $G$ over $K$, i.e.\ $\rho$ is an injective group homomorphism
and $V$ is a vector space over $K$ with $\dim_K V<\infty$. Define
$\fn{ed}_K(G)=\min\{\fn{trdeg}_K E:E$ is a faithful $G$-subfield
of $K(V)\}$. It is known that $\fn{ed}_K(G)$ is independent of the
choice of the faithful representation (see \cite[Theorem
3.1]{BR}).
\end{defn}

\begin{lemma}\label{l2.4}
Let $K$ be a field with $\fn{char}K=p>0$. Suppose that $\sigma\in
PGL_2(K)$ and $\fn{ord}(\sigma)$ in $PGL_2(K)$ is finite. Then
either $p\nmid \fn{ord}(\sigma)$ or $\fn{ord}(\sigma)=p$.
\end{lemma}

\begin{proof}
Note that the order of $\sigma$ in $PGL_2(K)$ is the same as that
in $PGL_2(\overline{K})$ where $\overline{K}$ is the algebraic
closure of $K$.

Choose a matrix $T\in GL_2(K)$ such that its image in $PGL_2(K)$
is $\sigma$. Find the Jordan canonical form of $T$ in
$GL_2(\overline{K})$. It is of the form
\[
\left(\begin{array}{cc}  a & 1 \\ 0 & a  \end{array}\right) \mbox{
\ or \ } \left(\begin{array}{cc}  a & 0 \\ 0 & b
\end{array}\right)
\]
where $a,b\in\overline{K}\backslash\{0\}$. It is not difficult to
see that the order of the image of the above matrix in
$PGL_2(\overline{K})$ is $p$ for the first case, and the order of
the image in $PGL_2(\overline{K})$ is relatively prime to $p$ for
the second case.
\end{proof}

\begin{lemma}\label{l2.5}
{\rm (\cite[Lemma 7.2]{BF})} Let $K$ be an arbitrary field and $G$
be a finite group. If $\fn{ed}_K(G)=1$, then $G$ can be embedded
into $PGL_2(K)$.
\end{lemma}

\begin{lemma}\label{l2.6}
Let $K$ be an arbitrary field, $\sigma\in PGL_2(K)$ be an element
of finite order. If $n=\fn{ord}(\sigma)$ and $\fn{char}K\nmid n$,
then $\zeta_n+\zeta^{-1}_n\in K$.
\end{lemma}

\begin{remark}\rm
The above lemma is a generalization of \cite[Lemma 7.7]{BF} where
it was required that $n$ is a prime number.
\end{remark}

\begin{proof}
Choose a matrix $T\in GL_2(K)$ such that its image in $PGL_2(K)$
is $\sigma$. Then the Jordan canonical form of $T$ is
\[
\left(\begin{array}{cc}  \lambda \zeta_n & 0 \\ 0 & \lambda
\end{array} \right)
\]
for some $\lambda$ in the algebraic closure of $K$.

Note that the rational canonical form of $T$ is
\[
\left(\begin{array}{cc}  \lambda \zeta_n & 0 \\ 0 & \lambda
\end{array} \right) \mbox{ \ or \ } \left(\begin{array}{cc}  0 & a
\\ 1 & b  \end{array} \right)
\]
where $\lambda,a,b\in K$ (according to whether the characteristic
polynomial is reducible over $K$ or irreducible over $K$). Then
first possibility will imply $\zeta_n\in K$; in particular
$\zeta_n+\zeta^{-1}_n\in K$. It remains to consider the second
possibility.

Compare the traces and the determinants of these two canonical
forms. We find that $b=\lambda\zeta_n+\lambda$ and
$-a=\lambda^2\zeta_n$. Thus
$\zeta_n+\zeta^{-1}_n+2=(\lambda\zeta_n+\lambda)^2(\lambda^2\zeta_n)^{-1}=b^2\cdot
(-a)^{-1}\in K$.
\end{proof}

\begin{lemma}\label{l2.7}
Let $p$ be a prime number and $K$ be a field with
$\fn{char}K=p>0$. For any positive integer $r$,
$\fn{ed}_K((\bm{Z}/p\bm{Z})^r)=1$ if and only if $[K:\bm{F}_p]\ge
r$.
\end{lemma}

\begin{proof}
Suppose that $[K:\bm{F}_p]\ge r$.

Let $(\bm{Z}/p\bm{Z})^r=\langle \sigma_i:\sigma^p_i=1,\
\sigma_i\sigma_j=\sigma_j\sigma_i$ for $1\le i\le r\rangle$.
Choose $\alpha_1,\alpha_2,\ldots,\alpha_r\in K$ so that
$\alpha_1,\ldots,\alpha_r$ are linearly independent over
$\bm{F}_p$. Consider the faithful representation
$\rho:(\bm{Z}/p\bm{Z})^r\to GL_2(K)$ defined by
\[
\rho(\sigma_i)=\left(\begin{array}{cc}  1 & \alpha_i \\ 0 & 1
\end{array}\right)
\]
for $1\le i \le r$. Thus we have a faithful
$(\bm{Z}/p\bm{Z})^r$-field $K(x,y)$ provided by this
representation. Define $t=\frac{x}{y}$. Then $K(t)$ is a faithful
$(\bm{Z}/p\bm{Z})^r$-subfield again.

For the other direction, suppose that
$\fn{ed}_K((\bm{Z}/p\bm{Z})^r)=1$.

If $[K:\bm{F}_p]=\infty$, there is nothing to prove. So consider
the case that $K$ is a finite field $\bm{F}_q$ where $q=p^n$ for
some integer $n$.

Since $\fn{ed}_K((\bm{Z}/p\bm{Z})^r)=1$, we may embed
$(\bm{Z}/p\bm{Z})^r$ into $PGL_2(K)=PGL_2(\bm{F}_q)$ by Lemma
\ref{l2.5}. Let $f:GL_2(\bm{F}_q)\to PGL_2(\bm{F}_q)$ be the
canonical projection. The group $f^{-1}((\bm{Z}/p\bm{Z})^r)$ is an
extension of $(\bm{Z}/p\bm{Z})^r$ by $\bm{F}^{\times}_q$. Since
$\gcd \{|\bm{F}^{\times}_q|,|(\bm{Z}/p\bm{Z})^r|\}=1$, we find
that the group extension splits by Schur-Zassenhaus's Theorem
\cite[p.235]{Su}. Hence $(\bm{Z}/p\bm{Z})^r$ can be embedded into
$GL_2(\bm{F}_q)$.

Since $|GL_2(\bm{F}_q)|=q(q^2-1)(q-1)$, it follows that $q$ is
divisible by $p^r$.
\end{proof}

\begin{theorem}\label{t2.8}
\r{(\cite[Corollary 4.16]{BF})} If $p$ is a prime number and $K$
is a field such that $\zeta_p\in K$, then
$\fn{ed}_K((\bm{Z}/p\bm{Z})^r)=r$. In particular, for a field $K$
with $\fn{char}K\ne 2$, $\fn{ed}_K((\bm{Z}/2\bm{Z})^r)=r$.
\end{theorem}

\setcounter{equation}{0}
\section{Theorems of Klein and Dickson}

Most material in this section may be found in \cite{Sp,Su}.

\begin{defn}\label{d3.1}\rm
A binary dihedral group of order $4n$ is defined by the
presentation $\langle\sigma,\tau: \sigma^n=\tau^2,\
\tau^{-1}\sigma\tau=\sigma^{-1}\rangle$. If $K$ is an
algebraically closed field and $\fn{char}K\nmid 2n$, this group
can be embedded in $SL_2(K)$ (see \cite[p.89]{Sp}) by defining
\[
\sigma=\left(\begin{array}{cc} \zeta_{2n} & 0 \\ 0 & \zeta^{-1}_{2n} \end{array}\right),\ \
\tau=\left(\begin{array}{cc} 0 & \sqrt{-1} \\ \sqrt{-1} & 0 \end{array}\right).
\]
\end{defn}

\begin{defn}\label{d3.2}\rm
A binary tetrahedral (resp.\ octahedral, icosahedral) group $G$ is a central extension
of $A_4$ (resp.\ $S_4$, $A_5$) by $\bm{Z}/2\bm{Z}$,
i.e.\ there is an element $\sigma\in G$ so that (i) $\fn{ord}(\sigma)=2$,
(ii) $\sigma$ belongs to the center of $G$,
and (iii) $G/\langle \sigma \rangle$ is isomorphic to $A_4$ (resp.\ $S_4$, $A_5$).
\end{defn}

If $K$ is an algebraically closed field and $\fn{char}K\ne 2$,
then both the binary tetrahedral group and the binary octahedral group can be embedded in $SL_2(K)$
(see \cite[p.91--92]{Sp}).
As an abstract group,
the binary tetrahedral group is isomorphic to $SL_2(\bm{F}_3)$ and the binary octahedral group
is isomorphic to the representation group of $S_4$ in which the transpositions correspond to
the elements of order 4 \cite[Theorem 6.17, p.404]{Su}.

If $K$ is an algebraically closed field and $\fn{char}K\nmid 10$,
then the binary icosahedral group can be embedded in $SL_2(K)$ (see \cite[p.93]{Sp}).
As an abstract group,
the binary icosahedral group is isomorphic to $SL_2(\bm{F}_5)$.

\begin{theorem}\label{t3.3}
{\rm (Klein \cite[p.89--93; Su, Theorem 6.17, p.404]{Sp})} Let $K$
be an algebraically closed field and $\fn{char}K=0$. If $G$ is a
finite group, then $G$ can be embedded in $SL_2(K)$ if and only if
$G$ is isomorphic to a cyclic group, a binary dihedral group, a
binary tetrahedral group, a binary octahedral group or a binary
icosahedral group.
\end{theorem}

\begin{defn}\label{d3.4}\rm
We will define the group $G(n,p^r)$ with the condition that
$p\nmid n$, $s\mid r$ where $s:=[\bm{F}_p(\zeta^2_n):\bm{F}_p]$.
We will adopt the following convention: Whenever we talk about
$G(n,p^r)$, it is assumed that the condition $p\nmid n$ and $s\mid
r$ with $s:=[\bm{F}_p(\zeta^2_n):\bm{F}_p]$ is satisfied
automatically.

We will define the group $G(n,p^r)$ first as a subgroup of
$SL_2(K)$ where $K$ is an algebraically closed field with
$\fn{char}K=p>0$. Then another definition of $G(n,p^r)$ as an
abstract group will be given in the form of generators and
relations (see (\ref{eq3.2})). Finally the group $G(n,p^r)$ will
be characterized as a subgroup of $SL_2(K)$ (where $K$ is an
algebraically closed field with $\fn{char}K=p>0$), which is a
semi-direct product of an elementary abelian $p$-group with a
cyclic group (see Lemma \ref{l3.5}).
\end{defn}

Now suppose that $K$ is an algebraically closed field with
$\fn{char}K=p>0$. Regard $K$ as a vector space over $\bm{F}_p$.
Since $\zeta_n \in K$, $K$ is also a vector space over
$\bm{F}_p(\zeta_n)$ (and therefore over $\bm{F}_p(\zeta_n^2)$).
Choose a vector subspace $V$ of $K$ over $\bm{F}_p(\zeta^2_n)$ so
that $[V:\bm{F}_p]=r$. (Note that
$r=[V:\bm{F}_p(\zeta^2_n)][\bm{F}_p(\zeta^2_n):\bm{F}_p]$.) Choose
a basis $\alpha_1,\ldots,\alpha_r$ of $V$ over $\bm{F}_p$. Define
$\sigma_1,\ldots,\sigma_r$, $\tau\in SL_2(K)$ by
\begin{equation}
\sigma_i=\left(\begin{array}{cc} 1 & \alpha_i \\ 0 & 1 \end{array}\right),\ \
\tau=\left(\begin{array}{cc} \zeta_n & a \\ 0 & \zeta^{-1}_n \end{array}\right)
\label{eq3.1}
\end{equation}
where $a$ is any element in $K$ if $n\ge 3$,
while $a=0$ if $n=1$ or 2.

Define $G(n,p^r)$ to be the subgroup of $SL_2(K)$ generated by
$\sigma_1,\ldots,\sigma_r$, $\tau$, i.e.
$G(n,p^r)=\langle\sigma_1,\ldots,\sigma_r,\ \tau\rangle$. Note
that $G(1,p^r)$ is an elementary abelian $p$-group and $G(2,p^r)$
is a direct product of an elementary abelian $p$-group with
$\bm{Z}/2\bm{Z}$.

Define $Q=\langle\sigma_1,\ldots,\sigma_r\rangle \subset G(n,p^r)$.
It is clear that $Q$ is a normal subgroup of $G(n,p^r)$ and $Q$ is an elementary abelian $p$-group.
A typical element in $Q$ is of the form

\[
\sigma=\left(\begin{array}{cc} 1 & v \\ 0 & 1 \end{array}\right)
\]
for some $v\in V$.
It is easy to verify that
\[
\tau\cdot\left(\begin{array}{cc} 1 & v \\ 0 & 1 \end{array}\right) \cdot \tau^{-1}
=\left(\begin{array}{cc} 1 & \zeta^2\cdot v \\ 0 & 1 \end{array}\right).
\]

The next step is to define $G(n,p^r)$ as an abstract group. Choose
a basis $\beta_1,\ldots,\beta_t$ of $V$ over $\bm{F}_p(\zeta^2_n)$
(thus $r=st$ where $s=[\bm{F}_p(\zeta^2_n):\bm{F}_p]$). Let
$f(X)=X^s-a_s X^{s-1}-a_{s-1}X^{s-2}-\cdots-a_1\in\bm{F}_p[T]$ be
the minimum polynomial of $\zeta^2_n$ over $\bm{F}_p$. (Note that
$f(X)$ is an irreducible factor of the cyclotomic polynomial
$\Phi_n(X)$ or $\Phi_{n/2}(X)$ over $\bm{F}_p$.) Define
$\beta_{ij}=\zeta_n^{2(j-1)}\beta_i$ where $1\le j\le s$. Then
$\beta_{ij}$ is a basis of $V$ over $\bm{F}_p$. It is not
difficult to show that $G(n,p^r)$ is generated by
\[
\left(\begin{array}{cc} 1 & \beta_{ij} \\ 0 & 1 \end{array}\right),\ \
\left(\begin{array}{cc} \zeta_n & a \\ 0 & \zeta^{-1}_n \end{array}\right)
\]
where $1\le i\le t$, $1\le j\le s$ and $r=st$. Moreover, the group
$G(n,p^r)$ may be defined by generators $\sigma_{ij}$ and $\tau$
(with $1\le i\le t$, $1\le j\le s$) and the relations are given by
\begin{eqnarray}
&& \sigma^p_{ij}=\tau^n=1,\ \ \sigma_{ij}\sigma_{kl}=\sigma_{kl}\sigma_{ij}, \nonumber \\
&& \tau\sigma_{ij}\tau^{-1}=\sigma_{i,j+1}\ \mbox{ \ for \ } 1\le i\le t,\ 1\le j\le s-1, \label{eq3.2} \\
&& \tau\sigma_{i,s}\tau^{-1}=\prod_{1\le j\le s} \sigma_{i,j}^{a_j}\ \mbox{ \ for \ }1\le i\le t. \nonumber
\end{eqnarray}

Thus, as an abstract group,
$G(n,p^r)$ is independent of the choice of $a$ in (\ref{eq3.1}).

\begin{lemma}\label{l3.5}
Let $K$ be an algebraically closed field with $\fn{char}K=p>0$.
Let $G$ be a finite subgroup of $SL_2(K)$ and $Q$ be a $p$-Sylow
subgroup of $G$. Assume that $Q$ is a normal subgroup of $G$ and
$Q$ is an elementary abelian group so that $G/Q$ is a cyclic
group. Then $G$ is conjugate to $G(n,p^r)$ for some integers $n$
and $r$.
\end{lemma}

\begin{proof}
Note that $G$ is a semi-direct product of $Q$ and $\langle\tau
\rangle$ ($\simeq G/Q$) because $\gcd\{|Q|,$ $|G/Q|\}=1$ and we
may apply Schur-Zassenhaus Theorem \cite[Theorem 8.10, p.235]{Su}.
Let $n=\fn{ord}(\tau)$.

Since $Q$ is an elementary abelian group, we may triangulate all
elements of $Q$ simultaneously. In other words, up to conjugation
in $SL_2(K)$, we may assume elements of $Q$ are of the form
\[
\left(\begin{array}{cc} 1 & v \\ 0 & 1 \end{array}\right)
\]
for some $v\in K$.
Write $Q=\langle\sigma_1,\ldots,\sigma_r\rangle$ ($\simeq (\bm{Z}/p\bm{Z})^r$) and
\[
\sigma_i=\left(\begin{array}{cc} 1 & \alpha_i \\ 0 & 1 \end{array}\right)
\]
for $\alpha_i\in K$.
Define $V=\bigoplus_{1\le i \le r}\bm{F}_p\cdot \alpha_i\subset K$.
Then, any $\sigma\in Q$ can be written as
\[
\sigma=\left(\begin{array}{cc} 1 & v \\ 0 & 1 \end{array}\right)
\]
for some $v\in V$.

Write
\[
\tau=\left(\begin{array}{cc} a & b \\ c & d \end{array}\right) \in SL_2(K).
\]

Since $\tau\sigma\tau^{-1}\in Q$ for any $\sigma\in Q$,
it follows that $c=0$.
Thus both $a$ and $d$ are primitive $n$-th root of unity.
Hence, without loss of generality,
we may write
\[
\tau=\left(\begin{array}{cc} \zeta_n & b \\ 0 & \zeta^{-1}_n \end{array}\right).
\]

Plugging in the relation $\tau\sigma\tau^{-1}\in Q$ for any $\sigma\in Q$,
we get
\[
\left(\begin{array}{cc} \zeta_n & b \\ 0 & \zeta^{-1}_n \end{array}\right)
\left(\begin{array}{cc} 1 & v \\ 0 & 1 \end{array}\right)
\left(\begin{array}{cc} \zeta_n & b \\ 0 & \zeta^{-1}_n \end{array}\right)^{-1}=
\left(\begin{array}{cc} 1 & \zeta^2_nv \\ 0 & 1 \end{array}\right)
\]
for any $v\in V$. We find that $\zeta^2_n\cdot V \subset V$, i.e.\
$V$ admits a structure of vector space over $\bm{F}_p(\zeta^2_n)$.
Hence $G$ is equal to the group $G(n,p^r)$ defined in (3.1).
\end{proof}

\begin{theorem}\label{t3.6}
{\rm (Dickson \cite[Theorem 6.17, p.404]{Su})}
Let $K$ be an algebraically closed field with $\fn{char}K=p>0$.
If $G$ is a finite group,
then $G$ can be embedded in $SL_2(K)$ if and only if $G$ is isomorphic to one of the following groups

Case I. When $p\nmid |G|$
\leftmargini=18mm
\begin{enumerate} \itemsep=-1pt
\item[\r{(i)}]
A cyclic group.
\item[\r{(ii)}]
A binary dihedral group of order $4n$.
\item[\r{(iii)}]
The binary tetrahedral group, i.e.\ $SL_2(\bm{F}_3)$.
\item[\r{(iv)}]
The binary octahedral group,
i.e.\ the representation group of $S_4$ in which the transpositions correspond to elements of order $4$.
\item[\r{(v)}]
The binary icosahedral group, i.e.\ $SL_2(\bm{F}_5)$.
\end{enumerate}

Case II. When $p\mid |G|$
\begin{enumerate} \itemsep=-1pt
\item[\r{(vi)}]
The group $G(n,p^r)$ with $p\nmid n$ and $s\mid r$ where $s=[\bm{F}_p(\zeta^2_n):\bm{F}_p]$.
\item[\r{(vii)}]
$p=2$ and $G=D_n$ with $n$ being an odd integer.
\item[\r{(viii)}]
$p=3$ and $G=SL_2(\bm{F}_5)$.
\item[\r{(ix)}]
$q$ is a power of $p$ and $G=SL_2(\bm{F}_q)$.
\item[\r{(x)}]
$p$ is odd, $q$ is a power of $p$ and
\[
G=\left\langle SL_2(\bm{F}_q),\
\left(\begin{array}{cc} \varepsilon & 0 \\ 0 & \varepsilon^{-1} \end{array}\right) \right\rangle
\]
where $\varepsilon\in K$ satisfies that
$\bm{F}_q(\varepsilon)=\bm{F}_{q^2}$ and
$\bm{F}^{\times}_q=\langle\varepsilon^2\rangle$.
\end{enumerate}
\end{theorem}

\begin{remark}\rm
In the statement of the above theorem, once we say the group in
(iii) is a finite subgroup of $SL_2(K)$, it is assumed tacitly
that $p=\fn{char}K\ne 2$ or 3. Note that, in (x), it is impossible
that $p=2$; otherwise, the condition
$\bm{F}_q(\varepsilon)=\bm{F}_{q^2}$ and
$\bm{F}^{\times}_q=\langle\varepsilon^2\rangle$ would lead to a
contradiction.
\end{remark}

\setcounter{equation}{0}
\section{Proof of Theorem \ref{t1.2}}

We will prove Theorem \ref{t1.2} in this section by using Theorem
\ref{t3.3} and Theorem \ref{t3.6}.

Let $K$ be an arbitrary field, $\overline{K}$ its algebraic
closure.

Suppose that $\fn{ed}_K(G)=1$. Then $G$ may be embedded in
$PGL_2(K)$ by Lemma \ref{l2.5}. Since $PGL_2(K)\subset
PGL_2(\overline{K})= PSL_2(\overline{K})$. We may regard $G$ as a
subgroup of $PSL_2(\overline{K})$.

Let $\pi:SL_2(\overline{K})\to PSL_2(\overline{K})$ be the natural
projection and $G'=\pi^{-1}(G)$. Then $G'$ is a finite subgroup of
$SL_2( \overline{K})$ and $G=\pi(G')$. The possible candidates for
$G'$ are prescribed in Theorem \ref{t3.3} and Theorem \ref{t3.6}.

Case 1. $\fn{char}K=0$.

Apply Theorem \ref{t3.3}.

If $G'$ is a cyclic group,
then $G$ is a cyclic group also.

If $G'$ is a binary dihedral group of order $4m$,
then $G$ is a dihedral group of order $2m$.

If $G'$ is a binary tetrahedral (resp.\ octahedral, icosahedral)
group, then $G$ is isomorphic to $A_4$ (resp.\ $S_4$, $A_5$).
Since $A_4 \subset S_4$ and $A_4 \subset A_5$, it follows that
$\fn{ed}_K(A_4)\le \fn{ed}_K(S_4)$, $\fn{ed}_K(A_4)\le
\fn{ed}_K(A_5)$. Because $\fn{ed}_K(A_4)\ge
\fn{ed}_K((\bm{Z}/2\bm{Z})^2)=2$ by Theorem \ref{t2.8}. It follows
that $\fn{ed}_K(G)\ge 2$.

Case 2. $\fn{char}K=p>0$ and $p\ne 2$.

Apply Theorem \ref{t3.6}.

When $p\nmid |G|$, apply similar arguments as in Case 1.

Suppose $p\mid |G|$.

Note that the order of $G'$ is even.

Suppose that $G'=G(n,p^r)$. Then $n$ is even. Thus $G=\pi(G')$ is
of the form $G(n/2,p^r)$ by (3.1) and Lemma \ref{l3.5}.

If $p=3$ and $G'=SL_2(\bm{F}_5)$, then $G$ is isomorphic to $A_5$.
But $\fn{ed}_K(A_5) \ge \fn{ed}_K(A_4) \ge 2$.

If $p$ is an odd prime number and $G'=SL_2(\bm{F}_q)$ or a group
containing $SL_2(\bm{F}_q)$, then $G$ contains $PSL_2(\bm{F}_q)$.
By \cite[Exercise 5(c), p.417]{Su}, $PSL_2(\bm{F}_q)$ contains
$(\bm{Z}/2\bm{Z})^2$. But $\fn{ed}_K((\bm{Z}/2\bm{Z})^2)=2$ by
Theorem \ref{t2.8}. Thus $\fn{ed}_K(G)\ge
\fn{ed}_K(PSL_2(\bm{F}_q))\ge 2$.

Case 3. $\fn{char}K=2$.

Apply Theorem \ref{t3.6} again.

Note that $G'\simeq G$ and thus $G$ may be regarded as a finite subgroup of $SL_2(\overline{K})$.

The possible candidates for $G$ are cyclic groups and groups in Case II of Theorem \ref{t3.6}.
The groups in Case II of Theorem \ref{t3.6} are
(vi) $G(n,2^r)$ with $n$ being an odd integer,
(vii) $D_n$ with $n$ being an odd integer,
and (ix) $SL_2(\bm{F}_q)$ where $q$ is a power of 2.
All these groups belong to the list of Theorem \ref{t1.2}.
\hfill$\square$

\setcounter{equation}{0}
\section{Proof of Theorem \ref{t1.3}}

We will prove Theorem \ref{t1.3} in this section.

\begin{lemma}\label{l5.1}
Let $K$ be any arbitrary field with $\fn{char}K\nmid n$.
If $n$ is an even integer and $\fn{ed}_K(\bm{Z}/n\bm{Z})=1$,
then $\zeta_n\in K$.
\end{lemma}

\begin{proof}
Step 1.
By Lemma \ref{l2.5} we may embed $G:=\bm{Z}/n\bm{Z}$ into $PGL_2(K)$.
Thus $PGL_2(K)$ has an element of order $n$.
By Lemma \ref{l2.6}, $\zeta_n+\zeta^{-1}_n\in K$.

Write $\eta=\zeta_n+\zeta^{-1}_n\in K$.
Define a matrix $T$ by
\[
T=\left(\begin{array}{cc} 0 & -1 \\ 1 & \eta \end{array}\right)\in GL_2(K).
\]

Note that the Jordan canonical form of $T$ in the algebraic closure of $K$ is
\[
\left(\begin{array}{cc} \zeta_n & 0 \\ 0 & \zeta^{-1}_n \end{array}\right).
\]

Thus the order of $T$ in $GL_2(K)$ is $n$,
while the order of the image of $T$ in $PGL_2(K)$ is $n/2$ because $n$ is even.

Define a faithful representation of $G=\langle\sigma\rangle$ to
$GL_2(K)$ by sending $\sigma$ to $T$. Let $K(x,y)$ be the
$G$-field associated with this representation so that $\sigma\cdot
x=\eta x+y$, $\sigma\cdot y=-x$ where $x$ and $y$ are
algebraically independent over $K$. (In fact, the action of
$\sigma$ on $x,y$ is given by the transpose inverse of the matrix
$T$.)

Step 2. Since $\fn{ed}_K(G)=1$, there is a faithful $G$-subfield
$E$ of $K(x,y)$ with $\fn{trdeg}_KE=1$. By L\"{u}roth's Theorem
$E$ may be written as $E=K(u)$ for some $u\in K(x,y)\backslash K$.

Write $u=g/f$ where $f,g\in K[x,y]$. We will find a generator
$w\in K(u)$, i.e.\ $K(u)=K(w)$, satisfying that $w=g_1/f_1$ where
$f_1,g_1\in K[x,y]$ and $\deg f_1\ne \deg g_1$ (note that $\deg
f_1$ and $\deg g_1$ denote the total degree of $f_1$ and $g_1$
with respect to $x$ and $y$).

Start from $u=g/f$.
If $\deg f\ne \deg g$, we are done.

Thus we may assume that $\deg f=\deg g$ in Step 3.

Step 3.
Note that $\sigma\in\fn{Aut}_K(K(u))\simeq PGL_2(K)$.
Hence $\sigma\cdot u=(au+b)/(cu+d)$ for some $a,b,c,d\in K$ with $ad-bc\ne 0$.

Substitute $u=g/f$ into the relation $\sigma\cdot u=(au+b)/(cu+d)$.
We get
\begin{equation}
(\sigma\cdot g)[cg+df]=(\sigma\cdot f)[ag+bf].  \label{eq5.1}
\end{equation}

Write $f=f_N+f_{N-1}+\cdots$, $g=g_N+g_{N-1}+\cdots$ where $N=\deg
f=\deg g$ and $f_i$, $g_i$ are homogeneous polynomials in $x$, $y$
of degree $i$. From (\ref{eq5.1}) we get
\begin{equation}
(\sigma\cdot g_N+\sigma\cdot g_{N-1}+\cdots)[(cg_N+df_N)+\cdots]
=(\sigma\cdot f_N+\sigma\cdot f_{N-1}+\cdots) [(ag_N+bf_N)+\cdots]  \label{eq5.2}
\end{equation}
since $\sigma$ is a linear map on $Kx+Ky$.

We claim that $(cg_N+df_N)(ag_N+bf_N)\ne 0$.

Otherwise, assume that $cg_N+df_N=0$.
Compare the degrees of both sides of (\ref{eq5.2}).
We find that $ag_N+bf_N=0$.
Since $g_N\cdot f_N\ne 0$,
it follows that $ad-bc=0$.
A contradiction.

Since $(\sigma\cdot g_N)(\sigma\cdot f_N)\cdot (cg_N+df_N)\cdot (ag_N+bf_N)\ne 0$,
the degrees of both sides of (\ref{eq5.2}) are $2N$.
Compare the leading terms of (\ref{eq5.2}).
We get
\begin{equation}
(\sigma\cdot g_N)[cg_N+df_N]=(\sigma\cdot f_N)[ag_N+bf_N].  \label{eq5.3}
\end{equation}

Write $\lambda=g_N/f_N$ and re-write (\ref{eq5.3}) as $\sigma\cdot\lambda=(a\lambda+b)/(c\lambda+d)$.
Note that the order of
\[
S:=\left(\begin{array}{cc} a & b \\ c & d \end{array}\right)\in PGL_2(K)
\]
is $n$ because $\sigma\cdot u=(au+b)/(cu+d)$ and $\sigma$ is faithful on $K(u)$.

We claim that $\lambda\in K\backslash\{0\}$.

Otherwise $\lambda$ is transcendental over $K$ and $\sigma$ is faithful on $K(\lambda)$
because the order of $S$ in $PGL_2(K)$ is $n$.

On the other hand write $f_N=\sum_{0\le i \le N} a_i x^{N-i} y^i$,
$g_N=\sum_{0\le i \le N} b_i x^{N-i} y^i$, $t=y/x$. Then
$\lambda=g_N/f_N=(\sum b_it^i)/(\sum a_i t^i) \in K(t)$.

Note that $\sigma\cdot t=\sigma\cdot y/\sigma\cdot x = -x/(\eta
x+y) =-1/(t+\eta)$. The order of the fractional linear
transformation $t\mapsto -1/(t+\eta)$ is $n/2$ because it is the
image of $T$ in $PGL_2(K)$. Hence $\sigma$ is not faithful on
$K(t)$. It follows that $\sigma$ is not faithful on $K(\lambda)$
because $K(\lambda)\subset K(t)$. A contradiction.

We conclude that $g_N/f_N=\lambda\in K\backslash \{0\}$.
Hence $u=g/f=(g_N+g_{N-1}+\cdots)/(f_N+f_{N-1}+\cdots)=\lambda+[(h_{N-1}+h_{N-2}+\cdots)/(f_N+f_{N-1}+\cdots)]$,
i.e.\ $u-\lambda=h/f$ with $f,h\in K[x,y]$ and $\deg f>\deg h$.
Clearly $K(u)=K(u-\lambda)$,
i.e.\ the goal of Step 2 is achieved.

Step 4.
In summary we find a faithful $G$-subfield $K(u)$ where $u=g/f$ and $\deg f \ne \deg g$.
Without loss of generality we may assume that $\deg f<\deg g$.

Since $\sigma\in \fn{Aut}(K(u))$,
$\sigma\cdot u=(au+b)/(cu+d)$ where $a,b,c,d\in K$ and $ad-bc\ne 0$.
We will show that $c=0$.

Write $g=g_N+g_{N-1}+\cdots$, $f=f_{N-1}+\cdots$ where $N=\deg g$
and $f_i$, $g_i$ are homogeneous polynomials in $x$, $y$.
Substitute it into $\sigma\cdot u=(au+b)/(cu+d)$ with $u=g/f$.
We get
\begin{equation}
(\sigma\cdot g)[cg_N+(c g_{N-1}+d f_{N-1})+\cdots]=
(\sigma\cdot f)[ag_N+(ag_{N-1}+b f_{N-1})+\cdots].  \label{eq5.4}
\end{equation}

If $c\ne 0$, the degree of the left-hand-side of (\ref{eq5.4}) is
$2N$ while that of the right-hand-side of (\ref{eq5.4}) is $\le
2N-1$. We conclude that $c=0$.

Thus we may write $\sigma\cdot u=\alpha\cdot u+\beta$ with $\alpha, \beta\in K$ and $\alpha\ne 0$.
It is easy to verify that $\sigma^i\cdot u=\alpha^i u+\beta(1+\alpha+\cdots+\alpha^{i-1})$ for any $i\ge 1$.
Since the order of $\sigma$ on $K(u)$ is $n$,
we find that $\alpha$ is a primitive $n$-th root of unity,
i.e.\ $\zeta_n\in K$.
\end{proof}

\begin{proof}[\indent Proof of Theorem \ref{t1.3}]

Let $G=\bm{Z}/n\bm{Z}=\langle\sigma\rangle$.

Suppose that $\fn{ed}_K(G)=1$.
If $\fn{char}K\nmid n$,
then $\zeta_n+\zeta^{-1}_n\in K$ by Lemma \ref{l2.5} and Lemma \ref{l2.6}.
Thus it is necessary that $\zeta_n+\zeta^{-1}_n\in K$ in the particular case when $n$ is odd.
When $n$ is even, apply Lemma \ref{l5.1}.
If $\fn{char}K=p>0$ and $p\mid n$,
we may apply Lemma \ref{l2.5} and Lemma \ref{l2.4} to conclude that $n=p$.

It remains to show that these necessary conditions are sufficient also.

When $\fn{char}K\nmid n$ and $n$ is odd, since
$\eta:=\zeta_n+\zeta^{-1}_n\in K$, we may define a faithful
representation of $G$ to $GL_2(K)$ as in the proof of Step 1 of
Lemma \ref{l5.1}. Let $K(x,y)$ be the same as in the proof of
Lemma \ref{l5.1}. Define $t=y/x$. Then $\sigma\cdot t=-1/t+\eta$.
The order of this fractional linear transformation is $n$ because
$n$ is an odd integer. Thus $K(t)$ is a faithful $G$-subfield of
$K(x,y)$. Hence $\fn{ed}_K(G)=1$.

When $\fn{char}K\nmid n$ and $n$ is even, since $\zeta_n\in K$,
we have a faithful one-dimensional representation of $G$.
Thus $\fn{ed}_K(G)=1$.

When $\fn{char}K=p>0$ and $p\mid n$,
since $n=p$, it is easy to see that
\[
\sigma\mapsto \left(\begin{array}{cc} 1 & a \\ 0 & 1 \end{array}\right)
\]
is a faithful representation for any $a\in K\backslash\{0\}$.
Consider the $G$ action on $K(x,y)$ given by $\sigma\cdot x=x-ay$,
$\sigma\cdot y=y$ with $\fn{trdeg}_K K(x,y)=2$. Define $t=x/y$. We
find $\sigma\cdot t=t-a$. Thus $K(t)$ is a faithful $G$-subfield.
It follows that $\fn{ed}_K(G)=1$.
\end{proof}

\setcounter{equation}{0}
\section{Proof of Theorem \ref{t1.4}}

Case 1. Assume that $\fn{char}K=0$.

If $\fn{ed}_K(D_n)=1$, by Lemma \ref{l2.5} we may embed $D_n$ into
$PGL_2(K)$. Thus $PGL_2(K)$ contains an element $\sigma$ of order
$n$. Apply Lemma \ref{l2.6} to get the necessary condition that
$\zeta_n+\zeta^{-1}_n\in K$.

We will show that $n$ is odd. Suppose to the contrary that $n$ is
even. Then $(\bm{Z}/2\bm{Z})^2\subset D_n$. Since
$\fn{ed}_K((\bm{Z}/2\bm{Z})^2)=2$ by Theorem \ref{t2.8}, we find
that $\fn{ed}_K(D_n)\ge 2$, which is a contradiction.

Now we consider the reverse direction. Assume that  $n$ is odd and
$\eta:=\zeta_n+\zeta^{-1}_n\in K$. Define matrices $T$ and $S$ in
$GL_2(K)$ as follows,
\[
T=\left(\begin{array}{cc} 0 & -1 \\ 1 & \eta
\end{array}\right), \ \ S=\left(\begin{array}{cc} 0 & 1 \\ 1 & 0
\end{array}\right).
\]

Let $\sigma$ and $\tau$ be generators of $D_n$ with relations
$\sigma^n=\tau^2=(\tau \sigma)^2=1$. Define a faithful
representation of $D_n$ into $GL_2(K)$ by sending $\sigma$ to $T$,
and $\tau$ to $S$. Then we have a faithful $D_n$-field $K(x,y)$
associated with this representation. Define $t=\frac{y}{x}$. Then
$K(t)$ is a faithful $D_n$-subfield because $n$ is odd. Thus
$\fn{ed}_K(D_n)=1$.

Case 2. Assume that $\fn{char}K=p > 0$ and $p\ne 2$.

Subcase 2.1.  Suppose $p\nmid n$.

The proof is the same as Case 1 because we may apply Lemma
\ref{l2.6}.

Subcase 2.2.  Suppose $p\mid n$.

If $\fn{ed}_K(D_n)=1$, then $PGL_2(K)$ contains an element
$\sigma$ of order $n$. Apply Lemma \ref{l2.4}. We find that $n=p$.

Conversely, let $D_p=\langle \sigma, \tau: \sigma^p=\tau^2=1,\
\tau\sigma\tau^{-1}=1\rangle$. Consider the faithful
representation $\rho:D_p\to GL_2(K)$ defined by
\[
\rho(\sigma)=\left(\begin{array}{cc}  1 & 1 \\ 0 & 1
\end{array}\right),\ \ \rho(\tau)=\left(\begin{array}{cc}  1 & 0
\\ 0 & -1  \end{array}\right).
\]

It is not difficult to show that $\fn{ed}_K(D_p)=1$.

Case 3. Assume that $\fn{char}K=2$.

Subcase 3.1. Suppose that $n$ is odd.

The situation is the same as in Case 1 or Subcase 2.1.

Subcase 3.2. Suppose that $n$ is even.

The situation is very similar to Subcase 2.2. If
$\fn{ed}_K(D_n)=1$, then $n=2$. Thus $D_n$ is isomorphic to
Klein's four group. Apply Lemma \ref{l2.7}. We find that $|K|\ge
4$.

Conversely, if $n=2$ and $|K|\ge 4$, choose $\alpha\in K\backslash
\{0,1\}$. Let $D_2=\langle \sigma,\tau: \sigma^2=\tau^2=1,\
\sigma\tau=\tau\sigma\rangle$. Define a faithful representation
$\rho:D_2\to GL_2(K)$ defined by
\[
\rho(\sigma)=\left(\begin{array}{cc}  1 & 1 \\ 0 & 1
\end{array}\right),\ \ \rho(\tau)=\left(\begin{array}{cc}  1 &
\alpha \\ 0 & 1  \end{array}\right).
\]

It is easy to show that $\fn{ed}_K(D_2)=1$.

\setcounter{equation}{0}
\section{Proof of Theorem \ref{t1.5}}

In this section $K$ is a field with $\fn{char}K=p>0$ and $G=G(n,p^r)$ where $p\nmid n$ and $s\mid r$
with $s:=[\bm{F}_p(\zeta^2_n):\bm{F}_p]$.

\begin{lemma}\label{l7.1}
Let $K$ be a field with $\fn{char}K=p>0$ and $p\ne 2$.
Let $G=G(n,p^r)$.
If $\fn{ed}_K(G)=1$,
then $\zeta_n\in K$ and $[K:\bm{F}_p]\ge r$.
\end{lemma}

\begin{proof}
Step 1.
By Lemma \ref{l2.5} we may embed $G$ into $PGL_2(K)$.
Since $PGL_2(K)\simeq \fn{Aut}_K(K(u))$ where $u$ is transcendental over $K$,
we may assume that $G$ acts faithfully on $K(u)$ by $K$-automorphisms.

Let $Q$ be the $p$-Sylow subgroup of $G$.
Then $Q$ is a normal subgroup of $G$ and $Q\simeq (\bm{Z}/p\bm{Z})^r$ (see Formula (\ref{eq3.2})).
Choose any $\sigma\in Q$, $\sigma\ne 1$.
Then $\sigma\cdot u=(au+b)/(cu+d)$ where $a,b,c,d\in K$ and $ad-bc\ne 0$.
We will find $w\in K(u)$ so that $K(u)=K(w)$ and $\sigma\cdot w=w+1$.

In fact, taking the rational canonical form of the matrix
\[
\left(\begin{array}{cc} a & b \\ c & d \end{array}\right)\in GL_2(K)
\]
amounts to finding another generator $w$ with $K(w)=K(u)$ and
$\sigma$ acting on $w$ according to this rational form. In other
words, without loss of generality, we may assume that the above
matrix is of its rational canonical form, i. e. it is of the form
\[
T=\left(\begin{array}{cc}  0 & -a \\ 1 & b \end{array}\right)\in GL_2(K).
\]

Thus $\sigma$ acts on $K(u)$ by $\sigma\cdot u=-a/(u+b)$ where
$a,b\in K$ and $a\ne 0$. Since $\fn{ord}(\sigma)=p$, the Jordan
canonical form of $T$ is
\[
\left(\begin{array}{cc} c & 1 \\ 0 & c \end{array}\right)
\]
for some $c\in \overline{K}\backslash \{0\}$.
It follows that
\begin{equation}
2c=b, \ \ c^2=a.  \label{eq7.1}
\end{equation}

Thus we find that $4a=b^2$.

Define $w=b/(2u+b)$. Since $\sigma\cdot u=-a/(u+b)$, it follows
that $\sigma\cdot w=w+1$.

Step 2. Once we know $\sigma\cdot w=w+1$, we can show that
$\lambda\cdot w=\alpha_\lambda w+\beta_\lambda$ for any
$\lambda\in Q$ where $\alpha_\lambda,\beta_\lambda\in K$ and
$\alpha_\lambda\ne 0$.

For, $\lambda\sigma=\sigma\lambda$ and $\lambda\cdot w =(\alpha
w+\beta)/(\gamma w+\delta)$ for some
$\alpha,\beta,\gamma,\delta\in K$ with
$\alpha\delta-\beta\gamma\ne 0$. From the relation
$\lambda\sigma(w)=\sigma\lambda(w)$, we get
\[
\left(\frac{(\alpha w+\beta)}{(\gamma w+\delta)}\right)+1
=\frac{[\alpha(w+1)+\beta]}{[\gamma(w+1)+\delta]} .
\]

It follows that $\gamma=0$.

In other words,
for any $\lambda \in Q$,
there exist $\alpha_\lambda, \beta_\lambda\in K$,
$\alpha_\lambda\ne 0$ so that $\lambda\cdot w=\alpha_\lambda w+\beta_\lambda$.

Since $\lambda^p=1$ for any $\lambda\in Q$,
we find that $\alpha_\lambda=1$.

In conclusion, for any $\lambda\in Q$, there is some
$\beta_\lambda\in K$ so that $\lambda\cdot w=w+\beta_\lambda$.
Moreover, it is easy to see the set $V=\{\beta_\lambda\in
K:\lambda\in Q\}$ is a vector space over $\bm{F}_p$ with
$[V:\bm{F}_p]=r$. Thus $[K:\bm{F}_p] \geq r$.

Step 3. Let $\tau\in G$ be an element of order $n$ so that
$G=\langle Q,\tau \rangle$ (see Formula (\ref{eq3.2})). Suppose
$\tau\cdot w=(Aw+B)/(Cw+D)$ for some $A,B,C,D\in K$ with $AD-BC\ne
0$. For any $\lambda\in Q$, since $\tau
\lambda\tau^{-1}=\lambda'\in Q$, we get
$\tau\lambda=\lambda'\tau$. Using the formulae $\lambda\cdot
w=w+\beta_\lambda$, $\lambda'\cdot w=w+\beta_{\lambda'}$, we get
\[
[(A+\beta_\lambda C)w+(B+\beta_\lambda D)] [Cw+(\beta_{\lambda'}C+D)]
=[Cw+D][Aw+(\beta_{\lambda'}A+B)].
\]

It follows that $C=0$. Thus we may write $\tau\cdot w=a\tau +b$
for some $a,b\in K$ and $a\ne 0$. From $\fn{ord}(\tau)=n$ and
$\tau^i\cdot w=a^i\tau+b(1+a+a^2+\cdots+a^{i-1})$ for $i\ge 1$, we
find that $a$ is a primitive $n$-th root of unity i.e.\
$\zeta_n\in K$.
\end{proof}

\begin{lemma}\label{l7.2}
Let $K$ be a field with $\fn{char}K=2$ and $G=G(n,2^r)$.
If $\fn{ed}_K(G)=1$, then $\zeta_n\in K$ and $[K:\bm{F}_p]\ge r$.
\end{lemma}

\begin{proof}
We use the same notation and the arguments as in Step 1 of the
proof of Lemma \ref{l7.1}. But Formula (\ref{eq7.1}) becomes
\begin{equation}
0=2c=b,\ \ c^2=a  \label{eq7.2}
\end{equation}

Thus $\sigma\cdot u=a/u$ for some $a\in K\backslash\{0\}$.
Note that it may happen that $a\in K^2$ or $a\notin K^2$,
because $c$ lies in the algebraic closure of $K$.

Case 1. $a\in K^2$.

Write $a=c^2$ with $c\in K\backslash \{0\}$.
Define $w=c/(u+c)$.
Then $\sigma\cdot w=w+1$.

Once we get $\sigma\cdot w=w+1$, Step 2 and Step 3 of the proof of
Lemma \ref{l7.1} work also. Hence $\zeta_n\in K$.

Since $(\bm{Z}/2\bm{Z})^r\simeq Q\subset G$,
it follows that $\fn{ed}_K((\bm{Z}/2\bm{Z})^r)=1$.
By Lemma \ref{l2.7} we find $[K:\bm{F}_2]\ge r$.

Case 2. $a\notin K^2$.

Define $F=K(\sqrt{a})$.
Then $\zeta_n\in F$ by Case 1.
It follows that $\zeta_n=\alpha+\beta\sqrt{a}$ for some $\alpha,\beta\in K$.
Thus $\zeta^2_n=\alpha^2+a\beta^2\in K$.
Since $n$ is odd by the definition of $G(n,2^r)$,
we find that $\zeta^2_n$ is also a primitive $n$-th root of unity.
Thus $\zeta_n\in K$.

The fact $[K:\bm{F}_2]\ge r$ may be proved as in Case 1.
\end{proof}

\begin{lemma}\label{l7.3}
Let $K$ be a field with $\fn{char}K=p>0$ and $p\ne 2$.
Let $G=G(n,p^r)$.
If $n$ is an even integer, then $\fn{ed}_K(G)\ge 2$.
\end{lemma}

\begin{proof}
Suppose not.
Assume that $\fn{ed}_K(G)=1$ and $n$ is even.
We will find a contradiction.

Step 1. By Lemma \ref{l7.1} we find that $\zeta_n\in K$ and
$[K:\bm{F}_p]\ge r$. We will find a faithful two-dimensional
representation of $G$.

Let $Q=\langle\sigma_1,\sigma_2,\ldots,\sigma_r\rangle \simeq (\bm{Z}/p\bm{Z})^r$ be the $p$-Sylow subgroup of $G$,
and $G=\langle Q,\tau \rangle$ where $\fn{ord}(\tau)=n$.

The field $K$ may be regarded as a vector space over
$\bm{F}_p(\zeta_n)$. Thus it is also a vector space over
$\bm{F}_p(\zeta^2_n)$. Write $s=[\bm{F}_p(\zeta^2_n):\bm{F}_p]$
and $r=st$. Since $[K:\bm{F}_p]\ge r$, it follows that
$[K:\bm{F}_p(\zeta^2_n)]\ge t$. Find a set of linearly independent
vectors $\alpha_1,\alpha_2,\ldots,\alpha_t$ over
$\bm{F}_p(\zeta^2_n)$. Define $V=\bigoplus_{1\le i\le t}
\bm{F}_p(\zeta_n^2)\cdot\alpha_i \subset K$. Choose a basis
$\beta_1,\ldots,\beta_r$ of $V$ over $\bm{F}_p$. Define a
representation $\rho:G\to GL_2(V)$ by
\[
\rho:\sigma_i\mapsto \left(\begin{array}{cc} 1 & -\beta_i \\ 0 & 1
\end{array}\right),\ \ \tau\mapsto \left(\begin{array}{cc} \zeta^{-1}_n
& 0 \\ 0 & \zeta_n \end{array}\right).
\]

Then $\rho$ is a faithful representation.
Let $K(x,y)$ be the $G$-field associated with this representation so that $\sigma_i\cdot x=x+\beta_i y$,
$\sigma_i\cdot y=y$, $\tau\cdot x=\zeta_n x$, $\tau\cdot y=\zeta^{-1}_n y$.

Since $\fn{ed}_K(G)=1$,
there is some element $u\in K(x,y)\backslash K$ so that $K(u)$ is a faithful $G$-subfield.

Step 2.
For any $\lambda\in G$,
$\lambda(u)=(a_\lambda u+b_\lambda)/(c_\lambda u+d_\lambda)$ where $a_\lambda,b_\lambda,c_\lambda,d_\lambda\in K$
and $a_\lambda d_\lambda-b_\lambda c_\lambda \ne 0$.
Write $u=g/f$ where $f,g\in K[x,y]$.
We will find $w\in K(u)$ so that $K(u)=K(w)$,
$\deg g_1\ne \deg f_1$ where $w=g_1/f_1$ with $f_1,g_1\in K[x,y]$.
The proof is the same as Step 3 of the proof of Lemma \ref{l5.1} by considering the action of $\tau$ on $K(u)$.
The details are omitted.

In conclusion,
without loss of generality,
we may write $u=g/f$ where $f,g\in K[x,y]$ and $\deg f<\deg g$.

Step 3.
Since $\lambda(u)=(a_\lambda u+b_\lambda)/(c_\lambda u+d_\lambda)$ for $a_\lambda,b_\lambda,c_\lambda,d_\lambda\in K$
with $a_\lambda d_\lambda-b_\lambda c_\lambda\ne 0$ for any $\lambda\in G$,
we apply the same arguments in Step 4 of the proof of Lemma \ref{l5.1} to conclude that $c_\lambda=0$
for any $\lambda\in G$.

In short, for any $\sigma\in Q$,
$\sigma\cdot u=u+b_\sigma$ for some $b_\sigma\in K$ while $\tau\cdot u=\zeta_n u+b$ for some $b\in K$.

Step 4.
Let $Q=\langle \sigma_{ij}: 1\le i\le t,\ 1\le j\le s \rangle$ where $r=st$ and $\sigma_{ij}$
are the elements defined in Formula (\ref{eq3.2}).
Recall that $f(X)=X^s-a_s X^{s-1}-\cdots -a_2X-a_1\in \bm{F}_p[X]$ is the minimum polynomial
of $\zeta^2_n$ over $\bm{F}_p$.

Write $\sigma_{ij}\cdot u=u+b_{ij}$ and $\tau\cdot u=\zeta_n u+b$.
We find that $\tau \sigma_{ij} \tau^{-1}\cdot u=u+\zeta_n^{-1}
b_{ij}$.

By Formula (\ref{eq3.2}),
$\tau\sigma_{ij}\tau^{-1}=\sigma_{i,j+1}$ if $1\le j\le s-1$, and
$\tau \sigma_{is}\tau^{-1}=\prod_j \sigma_{ij}^{a_j}$. Write
$\beta=b_{11}\ne 0$. We find that $\zeta^{-s}_n\beta=\sum_{1\le
j\le s} a_j \zeta^{-(j-1)}\beta$. We get that
$\zeta^{-s}_n-a_s\zeta^{-(s-1)}-a_{s-1}\zeta^{-(s-2)}-\cdots-a_2\zeta^{-1}-a_1=0$,
i.e.\ $f(\zeta^{-1}_n)=0$. But $f(X)$ is a factor of
$\Phi_{n/2}(X)$ (where $\Phi_{n/2}(X)$ is the $n/2$-th cyclotomic
polynomial) because $f(X)$ is the minimum polynomial of
$\zeta^2_n$. On the other hand $\zeta^{-1}_n$ is a primitive
$n$-th root of unity. We find that $f(X)$ divides
$\gcd\{\Phi_n(X),\Phi_{n/2}(X)\}$, which is impossible.
\end{proof}

\begin{proof}[\indent Proof of Theorem \ref{t1.5}]

Let $K$ be a field with $\fn{char}K=p>0$ and $G=G(n,p^r)$.

Suppose that $\fn{ed}_K(G)=1$.

If $p=2$,
then $n$ is odd by Definition \ref{d3.4}.

If $p\ne 2$ and $n$ is even, then it is impossible that
$\fn{ed}_K(G)=1$ by Lemma \ref{l7.3}. Hence $n$ is odd also.

The facts that $\zeta_n\in K$ and $[K:\bm{F}_p]\ge r$ follow from
Lemma \ref{l7.1} and Lemma \ref{l7.2}.

Conversely, assume that $n$ is odd,
$\zeta_n \in K$ and $[K:\bm{F}_p]\ge r$.
We will show that $\fn{ed}_K(G)=1$.

We will use the same notation and the same arguments in Step 1 of
the proof of Lemma \ref{l7.3}. In short, $K(x,y)$ is a faithful
$G$-field with $\sigma_i\cdot x=x+\beta_i y$, $\sigma_i\cdot y=y$,
$\tau\cdot x=\zeta_n x$, $\tau\cdot y =\zeta^{-1}_n y$ where $1\le
i \le r$ and $G=\langle \sigma_1,\ldots,\sigma_r,\ \tau\rangle$.

Define $t=x/y$.
Then $\sigma_i\cdot t=t+\beta_i$, $\tau\cdot t=\zeta^2_n t$.
Since $n$ is odd, $G$ acts faithfully on $K(t)$.
Hence $\fn{ed}_K(G)=1$.
\end{proof}

\begin{prop}\label{p7.4}
Let $K$ be a field with $\fn{char}K=2$.

\r{(1)} If $K\supset \bm{F}_4$, then $\fn{ed}_K(A_4)=\fn{ed}_K(A_5)=1$.

\r{(2)} If $K$ doesn't contain $\bm{F}_4$, then $\fn{ed}_K(A_4)=\fn{ed}_K(A_5)=2$.
\end{prop}

\begin{remark}\rm
Note that $\fn{ed}_K(A_3)=1$, because $\fn{ed}_K(S_3)=1$.
\end{remark}

\begin{proof}
(1) Note that $A_5\simeq SL_2(\bm{F}_4)$. If $K\supset\bm{F}_4$,
then we have a faithful $A_5$-field $K(x,y)$ provided by the
representation $A_5\rightarrow SL_2(\bm{F}_4) \subset GL_2(K)$.
Define $t=\frac{y}{x}$. Because $A_5$ is a simple group, it acts
faithfully on $K(t)$ . Hence $\fn{ed}_K(A_5)= 1$.

(2) First note that $A_4$ is isomorphic to $G(3,2^2)$. By Theorem
\ref{t1.5}, $\fn{ed}_K(G(3,2^2))=1$ if and only if $\zeta_3\in K$,
i.e.\ $K\supset \bm{F}_4$. In other words, if $K$ doesn't contain
$\bm{F}_4$, then $\fn{ed}_K(A_4)\ge 2$. On the other hand, we have
$\fn{ed}_K(A_5)\le \fn{ed}_K(S_5)\le 2$. It follows that $2\le
\fn{ed}_K(A_4) \le \fn{ed}_K(A_5)\le 2$.
\end{proof}

\setcounter{equation}{0}
\section{Proof of Theorem \ref{t1.6}}

In this section $K$ is a field with $\fn{char}K=2$ and $q$ is a power of 2.
Recall that $\zeta_{q+1}$ is a primitive $(q+1)$-th root of unity in the algebraic closure of $\bm{F}_2$.

\begin{lemma}\label{l8.1}
If $q=2^r$ for some positive integer $r$,
then $\bm{F}_2(\zeta_{q+1})=\bm{F}_{q^2}$ and $\bm{F}_2(\zeta_{q+1}+\zeta^{-1}_{q+1})=\bm{F}_q$.
\end{lemma}

\begin{proof}
Note that $\bm{F}_{q^2}=\bm{F}_2(\zeta_{q^2-1})$. Thus
$\zeta_{q+1}\in\bm{F}_2(\zeta_{q^2-1})=\bm{F}_{q^2}$. Hence
$\bm{F}_2(\zeta_{q+1})\subset \bm{F}_{q^2}$. To show that
$\bm{F}_2(\zeta_{q+1})=\bm{F}_{q^2}$, it suffices to show that
$\zeta_{q+1}$ doesn't belong to any proper subfield of
$\bm{F}_{q^2}$.

Since $[\bm{F}_{q^2}:\bm{F}_2]=2r$,
any proper subfield of $\bm{F}_{q^2}$ is of the form $\bm{F}_{2^m}$ where $m$ is a divisor of $2r$ and $m\ne 2r$.

Suppose that $2r=tm$ where $t\ge 2$ and $\zeta_{q+1}\in\bm{F}_{2^m}$.
Then $q+1$ divides $2^m-1$.
But $q+1=2^r+1$.
Since $m=2r/t\le r$, $2^r+1$ cannot be a divisor $2^m-1$.
Hence $\zeta_{q+1}\notin \bm{F}_{2^m}$ and $\bm{F}_2(\zeta_{q+1})=\bm{F}_{q^2}$.

Write $\eta=\zeta_{q+1}+\zeta^{-1}_{q+1}$. Then
$[\bm{F}_2(\zeta_{q+1}):\bm{F}_2(\eta)]\le 2$ because
$\zeta_{q+1}$ satisfies the equation $X^2-\eta X+1=0$. On the
other hand, note that $\zeta^{-1}_{q+1}=\zeta^q_{q+1}$. Hence
$\eta=\zeta_{q+1}+\zeta^{-1}_{q+1}=\zeta_{q+1}+\zeta^q_{q+1}=\zeta_{q+1}+\sigma(\zeta_{q+1})\in\bm{F}_q$
where $\sigma$ is the Frobenius automorphism of $\bm{F}_{q^2}$
over $\bm{F}_q$. Thus $\eta\in\bm{F}_q$. Hence
$\bm{F}_2(\eta)=\bm{F}_q$.
\end{proof}

\begin{proof}[\indent Proof of Theorem \ref{t1.6}]

Let $K$ be a field with $\fn{char}K=2$ and $G=SL_2(\bm{F}_q)$ where $q=2^r$ for some positive integer $r$.

If $r=1$, then $G=SL_2(\bm{F}_2)\simeq S_3$. It is known that
$\fn{ed}_K(S_3)=1$. Hence from now on we will assume that $r\ge
2$, and therefore $G$ is a simple group \cite[Theorem 9.9,
p.78]{Su}.

Suppose that $K\supset\bm{F}_q$. Then we have the trivial
representation of $G$ into $GL_2(K)$ by considering the inclusion
map $SL_2(\bm{F}_q)\subset GL_2(K)$. Let $K(x,y)$ be the
associated $G$-field with $\fn{trdeg}_K K(x,y)=2$ and $\sigma\cdot
x=ax+by$, $\sigma\cdot y =cx+dy$ if $\sigma\in G$ is of the form
\[
\sigma= \left(\begin{array}{cc} d & -b \\ -c & a
\end{array}\right) \in SL_2(\bm{F}_q).
\]

Define $t=x/y$.
If $\sigma$ acts on $x$, $y$ as above,
then $\sigma\cdot t=(at+b)/(ct+d)\in K(t)$.
Since $G$ is a simple group,
$G$ acts faithfully on $K(t)$.
Thus $\fn{ed}_K(G)=1$.

Conversely, assume that $\fn{ed}_K(G)=1$.

We claim that $G$ contains an element of order $q+1$. By Lemma
\ref{l8.1} $\bm{F}_2(\zeta_{q+1}+\zeta^{-1}_{q+1})=\bm{F}_q$, i.e.
$\eta:=\zeta_{q+1}+\zeta^{-1}_{q+1}\in \bm{F}_q$. Thus the
following matrix $T$ belongs to $SL_2(\bm{F}_q)$ where $T$ is
defined by
\[
T=\left(\begin{array}{cc} 0 & -1 \\ 1 & \eta \end{array}\right).
\]

The Jordan canonical form of $T$ is
\[
\left(\begin{array}{cc} \zeta_{q+1} & 0 \\ 0 & \zeta^{-1}_{q+1} \end{array}\right).
\]

Thus the order the $T$ is $q+1$.

Since $\fn{ed}_K(G)=1$, $G$ can be embedded in $PGL_2(K)$ by Lemma
\ref{l2.5}. Thus $PGL_2(K)$ contains an element of order $q+1$. By
Lemma \ref{l2.6}, we find $\zeta_{q+1}+\zeta^{-1}_{q+1} \in K$. By
Lemma \ref{l8.1},
$\bm{F}_q=\bm{F}_2(\zeta_{q+1}+\zeta^{-1}_{q+1})$. Hence
$\bm{F}_q\subset K$.
\end{proof}

\newpage

\end{document}